\documentclass[10pt,a4paper,twoside]{amsart}

\usepackage{amsfonts, amssymb, amsmath, amsthm}
\usepackage{mathrsfs} 
\usepackage{latexsym}
\usepackage{enumerate}
\usepackage{multicol}
\usepackage{verbatim}

\usepackage{url}
\usepackage{csquotes}

\usepackage[colorlinks=true]{hyperref}
\hypersetup{citecolor=blue, linkcolor=blue}
\usepackage{mathabx} 

\newtheorem{thm}{Theorem}
\newtheorem{prop}[thm]{Proposition}
\newtheorem{lem}[thm]{Lemma}
\newtheorem{cor}[thm]{Corollary}

\let\ve=\varepsilon

\def\Ocal{\mathcal{O}}

\def\1{1\!\!\!1}

\title{Accurate Computations of Euler Products over Primes in
  Arithmetic Progressions%
}
\author{Olivier Ramar\'e}
\address{CNRS / Aix Marseille Univ. / Centrale Marseille, I2M,
  Marseille, France}
\email{olivier.ramare@univ-amu.fr}
\subjclass[2010]{Primary 11Y60, Secundary 11N13, 05A}
\keywords{Euler products}

\begin{document}



\maketitle

\begin{abstract}
  \texttt{File \jobname.tex} 
  This note provides truncated formulae with  explicit error terms to
  compute Euler products over primes in arithmetic progressions of
  rational fractions. It further provides such a formula for the
  product of terms of the shape $F(1/p, 1/p^s)$ when $F$ is a
  two-variable polynomial with coefficients in $\mathbb{C}$ satisfying
  some restrictive conditions.
\end{abstract}


\section{Introduction and results}

Our primary concern in this note is to evaluate Euler products of the
shape
\begin{equation*}
  \prod_{p\equiv a[q]}\biggl(1-\frac{1}{p^s}\biggr)
\end{equation*}
when $s$ is a \emph{complex} parameter satisfying $\Re s>1$.
Such computations have attracted some attention as these values occur
when $s$ is a real number as densities in number theory.
D.~Shanks in \cite{Shanks*60} (resp. \cite{Shanks*61},
resp. \cite{Shanks*67}) has already computed accurately an Euler
product over primes congruent to~1 modulo~8 (resp. to~1 modulo~4,
resp. 1 modulo~8). His method has been extended by S.~Ettahri, L.~Surel
and the present author in~\cite{Ettahri-Ramare-Surel*19} in an
algorithm that converges very fast (double exponential convergence)
but this extension
covers only some special values for the residue class $a$, or some
special bundle of them; it is further limited to real values of~$s$. 

We 
use logarithms, and since the logarithm of a product is not a priori the sum of
the logarithms, we need to clarify things before embarking in this
project.
First, the $\log$-function corresponds in this paper always to what is
called \emph{the principal branch of the logarithm}. We recognize it
because its argument vanishes when we restrict it to the real line and
we consider it undefined on the non-positive real numbers. The second point
is contained 
in the next elementary proposition.
\begin{prop}
  \label{mylog}
  Associate to each prime $p$ a complex number $a_p$ such that
  $|a_p|<p$ and $a_p\ll_\ve p^\ve$ for every $\ve>0$. We consider the
  Euler product defined when $\Re s>1$ by:
  \begin{equation}
    \label{defD}
    D(s)=\prod_{p\ge2}\biggl(1-\frac{a_p}{p^s}\biggr)^{-1}.
  \end{equation}
  In this same domain we have
  \begin{equation}
    \label{deflogD}
    \log D(s)=\sum_{p\ge2}\sum_{k\ge1}\frac{a_p^k}{kp^{ks}}.
  \end{equation}
\end{prop}
This is simply because, by using the expansion of the principal branch
of the logarithm in Taylor series, namely $-\log(1-z)=\sum_{k\ge1}z^k/k$ valid
for any complex $z$ inside the unit circle, we find that
$C(s)=-\sum_{p\ge2}\log(1-a_p/p^s)$ verifies $\exp C(s)=D(s)$, so that
$C(s)$ is indeed a candidate for $\log D(s)$. The second remark is
that $D(s)$ approaches 1 when $\Re s$ goes to infinity while our
choice for $\log D(s)$ indeed approaches $0$ and no other multiple of
$2i\pi$. These two remarks are enough justification of this proposition.

We assume here that the values of the Dirichlet $L$-series
$L(s,\chi)$ may be computed with arbitrary precision when $\Re
s>1$. Our aim is thus to reduce our computations to these ones.
Here is an identity to do so.
\begin{thm}
  \label{identity}
  Let $a$ be prime to the modulus $q\ge1$ and let $\widehat{G}_q$ be the
  group of Dirichlet characters modulo $q$. We have
  \begin{equation*}
    -\mkern-10mu\sum_{\substack{p\equiv a[q],\\ p\ge P}}\log(1-1/p^s)
    =\sum_{\ell\ge1}\frac{-1}{\ell\varphi(q)}\sum_{d|\ell}\mu(d)
    \sum_{\chi\in \widehat{G}_q}\overline{\chi}(a)
    \log L_P(\ell s,\chi^d)
  \end{equation*}
  where
  \begin{equation}
  \label{defLP}
  L_P(s,\chi)=\prod_{p\ge P}(1-\chi(p)/p^s)^{-1}.
\end{equation}
\end{thm}
If finding this identity has not been immediate, checking it is only a
matter of calculations that we reproduce in
Section~\ref{checkthm}. A partial identity of this sort has already
been used by K.~Williams in \cite{Williams*74} and more recently by
A.~Languasco and A.~Zaccagnini in \cite{Languasco-Zaccagnini*10}.
It is worth noticing that, with our conventions, we have the obvious
\begin{equation*}
  \log L_P( s,\chi)
  =\log L( s,\chi)-\sum_{p<P}\log(1-\chi(p)/p^s).
\end{equation*}
This leads to the next immediate corollary.
\begin{cor}
  \label{identityshort}
  Let $a$ be prime to the modulus $q\ge1$ and let $\widehat{G}_q$ be the
  group of Dirichlet characters modulo $q$.
  Let further two integer parameters $P\ge2$ and $L\ge2$ be chosen. We have
  \begin{equation*}
    \prod_{\substack{p\ge P,\\ p\equiv
        a[q]}}\biggl(1-\frac{1}{p^s}\biggr)
    =
    \exp\biggl(Y_P(s;q,a| L)+\Ocal^*\biggl(\frac{1}{P^{L\Re s}}\biggr)\biggr)
    .
  \end{equation*}
  where
  \begin{equation}
    \label{defY}
    Y_P(s;q,a| L)
    =
    \sum_{\ell\le L}\frac{1}{\ell}\sum_{d|\ell}\mu(d)
    \sum_{\chi\in \widehat{G}_q}\frac{\overline{\chi}(a)}{\varphi(q)}
    \log L_P(\ell s,\chi^d)
  \end{equation}
\end{cor}

\subsubsection*{Extension to one variable rational fractions}
Once we have such an approximation, we can reuse the machinery
of~\cite{Ettahri-Ramare-Surel*19} to reach Euler products of the shape
\begin{equation*}
  \prod_{\substack{p\ge P,\\ p\equiv a[q]}}(1+R(p^s))
\end{equation*}
where $R$ is a rational fraction. 
\begin{thm}
\label{thm4}
  Let $F$ and $G$ be two polynomials of $\mathbb{C}[t]$. We assume
  that $G(0)=1$ and that $F(0)=F'(0)=0$.
  Let $\beta\ge2$ be larger than the inverse of the roots of $G$ and
  of $G-F$.
Let $P\ge 2\beta$ be an integer parameter.  Then, for any integer parameter $L\ge2$, we have
\begin{equation*}
  \prod_{\substack{p\ge P,\\ p\equiv a[q]}}
  \biggl(1-\frac{F(1/p)}{G(1/p)}\biggr)
  =
  \exp\biggl(\sum_{2\le j\le J}\bigl(b_{G-F}(j)-b_G(j)\bigr) Y_P(j;q,a|L)
  + I\biggr)
\end{equation*}
where the integers $b_{G-F}(j)$ and $b_G(j)$ are defined 
in Lemma~\ref{Wittpoly},
\begin{equation*}
    |I|\le 
    8\max(\deg (G-F),\deg G)
    \beta^2(\beta/P)^{2L}
  \end{equation*}
  and $Y(s;q,a|L)$ is defined by~\eqref{defY}.
\end{thm}
We obtained in \cite{Ettahri-Ramare-Surel*19} an approximation that is
much better but only valid for rational fractions with real
coefficients and some residue classes. 

One can write a similar theorem for the Euler product
\begin{equation*}
  \prod_{\substack{p\ge P,\\ p\in\mathcal{A}}}
  \biggl(1-\frac{F(1/p^s)}{G(1/p^s)}\biggr).
\end{equation*}

\subsubsection*{Extension to two variables rational fractions}
The general form of Euler products that one has
to treat in practice are of the shape
\begin{equation*}
  \prod_{\substack{p\ge P,\\ p\equiv a[q]}}(1+R(p, p^s))
\end{equation*}
where $R$ is a rational fraction of two variables.
When $s$ takes a specific rational value, typically $2$, $3/2$ or
$4/3$, this question reduces to the above one though each values of
$s$ requires a new rational fraction; this covers most of the cases
when we have to compute a single special constant.
In the general case however, for
instance when $s=2+i$, such a trick fails.
The theoretical
understanding of this situation is also limited even for
$q=1$. For instance, if the case of a rational fraction of one
variable is covered by the theorem of T.~Esterman in \cite{Esterman*28} and
extended by G.~Dalhquist in
\cite{Dahlquist*52}, no such result is known in the general situation.  
This question has been addressed in the context of enumerative algebra,
for instance by M.~du Sautoy and F.~Gr\"unewald in
\cite{duSautoy-Grunewald*02}. The lecture notes
\cite{duSautoy-Woodward*08} by M.~du Sautoy and L.~Woodward contains
material in this direction.There are several continuations of
Esterman's work; for instance, one may
 consider Euler products of the shape $R(p^{s_1},
p^{s_2})$ (with the hope of being able to specify $s_1$), see for instance \cite{Delabarre*13} by L.~Delabarre, but
these results do not apply to our case. 

We are able to handle some rational fractions by reducing them to the
case treated in the next theorem.
\begin{thm}
  \label{thm5}
  Let $s$ be a complex
  number with $\Re s=\sigma>1$. Let $(a_\ell)_{\ell\le k}$ be a
  sequence of complex numbers and $(u_\ell)_{\ell\le k}$ and
  $(v_\ell)_{\ell\le k}$ be two sequences of real numbers. We assume
  that $u_\ell\sigma+v_\ell>0$ and we define $A=\max(1,\max(|a_\ell|))$. Let $q$ be a
  modulus, $a$~be an invertible residue class modulo~$q$ and $P\ge2kA$
  and $L\ge k$
  be two integer parameters. We have
  \begin{equation*}
    \prod_{\substack{p\ge P,\\ p\equiv a[q]}}\biggl(
    1-\sum_{1\le \ell\le
      k}\frac{a_\ell}{p^{u_\ell s+v_\ell}}
    \biggr)
    =\exp{-(Z+I)}
  \end{equation*}
  where
  \begin{equation}
    \label{defZ}
    Z = \mkern-25mu\sum_{\substack{m_1,\ldots,m_k\ge 0,\\ 1\le m_1+\ldots+m_k\le
      L}}\mkern-30mu M(m_1,\ldots,m_k)\sum_{f\le
      F}
    \frac{\kappa_f(\prod_{\ell\le k}a_\ell^{m_\ell})}{f}
    Y_P\Bigl(\sum_{\ell\le k}m_\ell (u_\ell
    s+v_\ell);q,a|L\Bigr)
  \end{equation}
  where $M(m_1,m_2,\ldots,m_k)$ is defined at~\eqref{defM}, $\kappa_f$
  is defined at~\eqref{defkappa}, $Y(s;q,a|L)$ is defined
  by~\eqref{defY} and finally where
  \begin{equation}
    \label{eq:7}
    |I|\le \frac{2^k\cdot A^L}{k!P^L}\biggl(
    (L+k)^k
    +1+\log L+\frac{3kA}{L}
    \biggr).
  \end{equation}
\end{thm}
Hence this theorem provides us with an exponentially decreasing error term.
More complicated terms may be handled through this theorem by writing
\begin{align*}
  1+\frac{F(p,p^s)}{G(p,p^s)}
  &=
  \frac{(F+G)(p,p^s)}{p^{As+B}}\frac{p^{As+B}}{G(p,p^s)}
  \\&=
  \biggl(1+ \frac{(F+G)(p,p^s)-p^{As+B}}{p^{As+B}}\biggr)
  \biggl(1+\frac{G(p,p^s)-p^{As+B}}{p^{As+B}}\biggr)^{-1}.
\end{align*}
This would function when $G$ has a clearly dominant monomial. It typically
works for $G(p,p^s)=p^{2s}(p^2+1)$ but fails for $G(p,p^s)=p^{2s}(p+1)$.
Our most important additional tool, namely Lemma~\ref{MW}, may be used
to obtain results on analytic continuation, but since we use
logarithms elsewhere, the general effect is unclear. We however
provide the next example:
\begin{equation}
\label{inidiff}
    D(s)=\prod_{p\ge2}\biggl(1+\frac{1}{p^s}-\frac{1}{p^{2s-1}}\biggr).
\end{equation}
Lemma~\ref{MW} gives us the decomposition
\begin{equation*}
  D(s)=\prod_{\substack{m_1,m_2\ge0,\\ m_1+m_2\ge1}}
  \prod_{p\ge2}\biggl(1-\frac{(-1)^{m_1}}{p^{(m_1+2m_2)s-m_2}}\biggr)^{M(m_1,m_2)}.
\end{equation*}
We check that $M(1,0)=M(0,1)=1$ and that $M(m,0)=M(0,m)=0$ when
$m\ge2$, whence
\begin{equation}
  \label{eq:5}
  D(s)=\zeta(2s-1)\frac{\zeta(2s)}{\zeta(s)}
  \prod_{\substack{m_1,m_2\ge1}}
  \prod_{p\ge2}\biggl(1-\frac{(-1)^{m_1}}{p^{(m_1+2m_2)s-m_2}}\biggr)^{M(m_1,m_2)}.
\end{equation}
This writing offers an analytic continuation of $D(s)$ to the domain
defined by~$\Re s>1/2$. This analysis can be extended to
\begin{equation*}
  \prod_{p\ge2}\biggl(1-\frac{C_1}{p^s}-\frac{C_2}{p^{2s-1}}\biggr)
\end{equation*}
when $C_1$ and $C_2$ are \emph{integers}. In general, Lemma~\ref{MW}
transfers to problem to the analytic continuation of
$\prod_p(1-c/p^s)$ for some $c$ but even the case $c=\sqrt{2}$ is difficult.



  

  

\section{Proof of Theorem~\ref{identity} and its Corollary}
\label{checkthm}

\begin{proof}[Proof of Theorem~\ref{identity}]
We have to simplify the expression
\begin{equation}
  \label{eq:1}
  S=\sum_{\ell\ge1}\frac{1}{\ell\varphi(q)}\sum_{d|\ell}\mu(d)
    \sum_{\chi\in \widehat{G}_q}\overline{\chi}(a)
    \sum_{p\ge P}\sum_{k\ge1}\frac{\chi(p)^{dk}}{kp^{k\ell s}}.
  \end{equation}
  We readily check that, when $h\ge1$ and $p$ are
  fixed, we have
  \begin{align*}
    \sum_{k\ell=h}
    \sum_{d|\ell}\mu(d)
    \sum_{\substack{\chi\in\widehat{G}_q}}\overline{\chi}(a)\chi(p)^{dk}
    &=
      \sum_{k|h}
    \sum_{dk|h}\mu(d)
    \sum_{\substack{\chi\in\widehat{G}_q}}\overline{\chi}(a)\chi(p)^{dk}
    \\&=
      \sum_{g|h}
    \sum_{d|g}\mu(d)
    \sum_{\substack{\chi\in\widehat{G}_q}}\overline{\chi}(a)\chi(p)^{g}
    \\&=
    \sum_{\substack{\chi\in\widehat{G}_q}}\overline{\chi}(a)\chi(p)
    =\varphi(q)\1_{p\equiv a[q]}
  \end{align*}
  and the theorem follows directly.
\end{proof}

\begin{proof}[Proof of Corollary~\ref{identityshort}]
  A moment thought discloses that 
  \begin{equation*}
    |\log L_P(s,\chi)|\le \log \zeta_P(\sigma)
  \end{equation*}
  where $\sigma=\Re s$. We have furthermore
  \begin{equation*}
    \log \zeta_P(\sigma)
    \le \sum_{\substack{n\ge P}} \frac{1}{n^\sigma}
    \le \int_{P}^\infty\frac{dt}{t^{\sigma}}
    =\frac{1}{(\sigma-1)P^{\sigma-1}}.
  \end{equation*}
  by our assumptions. We next check that
  \begin{equation*}
    \biggl|
    \sum_{\ell> L}\frac{1}{\ell\varphi(q)}\sum_{d|\ell}\mu(d)
    \sum_{\chi\in \widehat{G}_q}\overline{\chi}(a)
    \log L_P(\ell s,\chi^d)
    \biggr|
    \le
    \sum_{\ell> L}\frac{2^{\omega(\ell)}}{\ell}
    \frac{P}{(\ell\sigma-1) P^{\ell\sigma}}.
  \end{equation*}
  Here $\omega(\ell)$ denotes the number of prime factors of $\ell$
  (without multiplicity).
  We use the simplistic bounds $2^{\omega(\ell)}\le \ell$ and
  $\ell\sigma-1\ge 2$. This yields the upper bound
  $\frac{P}{2P^{L\sigma}(P^{\sigma}-1)}$ which is no more than
  $1/P^{L\sigma}$. We finally recall that $e^x-1\le \frac{8}{7}x$ when
  $x\in[0,1/4]$ as the function $(e^x-1)/x$ is non-decreasing (its expansion
  in power series has non-negative coefficients).
\end{proof}

\section{Proof of Theorem~\ref{thm4}}
We first need to extend \cite[Lemma 16]{Ettahri-Ramare-Surel*19} to
cover the case of polynomials with complex coefficients.
The ancestor of this Lemma is \cite[Lemma 1]{Moree*00}.
\begin{lem}
\label{Wittpoly}
Let $H(t) = 1+a_1 t+\ldots+a_{\delta}t^{\delta} \in \mathbb{C}[t]$ be a 
polynomial of degree~$\delta$. 
Let $\alpha_{1},\ldots,\alpha_{\delta}$ be the inverses of its
roots. Put $s_{H}(k) =\alpha_{1}^{k}+\ldots+\alpha_{\delta}^{k}$. The
$s_{H}(k)$ satisfy the Newton-Girard recursion
\begin{equation}
\label{recursion}
s_{F}(k)+a_1s_F(k-1)+\ldots+a_{k-1}s_{F}(1)+ka_{k}=0,   
\end{equation}
where we have defined $a_{\delta+1} =a_{\delta+2}=\ldots=0$. We define
\begin{equation}
    \label{bfk}
b_{H}(k)=\frac{1}{k}\sum_{d|k}\mu({k}/{d})s_{H}(d).
\end{equation}
\end{lem}

\begin{lem}
  \label{step1}
  Let $F$ and $G$ be two polynomials of $\mathbb{C}[t]$. We assume
  that $G(0)=1$ and that $F(0)=0$.
  Let $\beta\ge1$ be larger than the inverse of the roots of $G$ and
  of $G-F$. When $z$
  is a complex number such that $|z|<\beta$ and $|1-(F/G)(z)|<1$. We have
  \begin{equation}
    \label{Fhatb}
    \log \biggl(1-\frac{F(z)}{G(z)}\biggr)
    =\sum_{j\ge1}\bigl(b_{G-F}(j)-b_G(j)\bigr)\log(1-z^j).
  \end{equation}
\end{lem}
%
\begin{proof}
  We adapt the proof of \cite[Lemma 1]{Moree*00}. We write
  $(G-F)(t)=\prod_{i}(1-\alpha_it)$. We have
  \begin{equation*}
    \frac{(G-F)'(t)}{(G-F)(t)}=\sum_i\frac{\alpha_i t}{1-\alpha_it}
    =\sum_{k\ge1}s_{G-F}(k)t^{k-1}.
  \end{equation*}
  This series is absolutely convergent in any disc $|t|\le b<1/\beta$
  where $\beta=\max_j(1/|\alpha_j|)$. 
  We may also decompose $(G-F)'(t)/(G-F)(t)$ in Lambert series as
  \begin{equation*}
    \frac{(G-F)'(t)}{(G-F)(t)}=\sum_{j\ge1}b_{G-F}(j)\frac{jt^{j-1}}{1-t^j}
  \end{equation*}
  as some series shuffling in any disc of radius $b<\min(1,1/\beta)$ shows.
  The comparison of the coefficients justify the formula~\eqref{bfk}.
  We may do the same for $G$ instead of $G-F$ (or use the case $F=0$).
  We find that
  \begin{equation*}
    \frac{G'-F'}{G-F}-\frac{G'}{G}
    =\frac{-(F'G-FG')}{G(G-F)}
    =\frac{-(F'G-FG')}{G^2}\sum_{k\ge0}\biggl(\frac{F}{G}\biggr)^k.
  \end{equation*}
  and by formal integration, this gives us the identity
  \begin{equation*}
    -\sum_{k\ge1}\frac{(F/G)(t)^k}{k}
    =-\sum_{j\ge1}\bigl(b_{G-F}(j)-b_G(j)\bigr)\log(1-t^j).
  \end{equation*}
  This readily extends into a equality between analytic function in
  the domain where $|(F/G)(z)-1|<1$ and $|z|<\beta$. The lemma follows readily.
\end{proof}

Here is now \cite[Lemma 17]{Ettahri-Ramare-Surel*19}, though for
polynomials with complex coefficients.
\begin{lem}\label{apriorimaj}
  We use the hypotheses and notation of Lemma~\ref{Wittpoly}. Let $\beta\ge2$ be larger than the inverse of the modulus of all the roots of $H(t)$. We have
  \begin{equation*}
      |b_H(k)|\le2\deg H \cdot \beta^k/k.
  \end{equation*}
\end{lem}
And we finally recall \cite[Lemma 18]{Ettahri-Ramare-Surel*19} that yields
an easy upper estimates for the inverse of the modulus of all the roots of $F(t)$ in terms of its coefficients.
\begin{lem}
\label{easybeta}
Let $H(X)=1+a_1X+\ldots+a_\delta X^\delta$ be a polynomial of
degree~$\delta$. Let $\rho$ be one of its roots. Then either $|\rho|\ge 1$ or
$1/|\rho|\le |a_1|+|a_2|+\ldots+|a_\delta|$.
\end{lem}

\begin{proof}[Proof of Theorem~\ref{thm4}]
The proof requires several steps. We start from Lemma~\ref{step1},
i.e. from the identity 
\begin{equation}
  \label{formal-FG}
  \log \biggl(1-\frac{F(z)}{G(z)}\biggr)
  =\sum_{j\ge2}\bigl(b_{G-F}(j)-b_G(j)\bigr)\log(1-z^j),
\end{equation}
in the domain $|z|<\beta$ and $|1-(F/G)(z)|<1$.
The fact that the term with $j=1$ vanishes comes from our assumption
that $F(0)=F'(0)=0$.
To control the rate of convergence, we
notice that 
By Lemma~\ref{apriorimaj}, we know that
$
    |b_{G-F}(j)-b_G(j)|\le 4\max(\deg (G-F),\deg G)\beta^j/j
    $.
Therefore, for any bound $J$, we have
\begin{equation}
   \label{tailJ}
    \sum_{j\ge J+1}|t^j||b_{G-F}(j)-b_G(j)|
    \le
    4\max(\deg (G-F),\deg G)\frac{|t\beta|^{J+1}}{(1-|t\beta|)(J+1)},
\end{equation}
as soon as $|t|<1/\beta$.    Furthermore, we deduce that $|\log(1-z)/z|\le \log(1-1/2)/(1/2)\le
    3/2$
    when $|z|\le 1/2$ by looking at the Taylor expansion.
 We thus have
\begin{equation}
  \label{true-FG}
  \log \biggl(1-\frac{F(z)}{G(z)}\biggr)
  =\sum_{2\le j\le J}\bigl(b_{G-F}(j)-b_G(j)\bigr)\log(1-z^j)
  +I_1
\end{equation}
where
$|I_1|\le 6\max(\deg (G-F),\deg G)|z\beta|^{J+1}/(1-|z\beta|)$.
Now that we have the expansion~\eqref{true-FG} at our disposal for each prime $p$, we may combine them. We readily get
\begin{equation*}
  \sum_{\substack{p\ge P,\\ p\equiv a[q]}}
  \log\biggl(1-\frac{F(1/p)}{G(1/p)}\biggr)
  =
  \sum_{2\le j\le J}\bigl(b_{G-F}(j)-b_G(j)\bigr)\sum_{\substack{p\ge
      P,\\ p\equiv a[q]}}
  \log(1-1/p^j)
    + I_2,
\end{equation*}
where $I_2$ satisfies
\begin{align*}
    | I_2|
    &\le
    6\max(\deg (G-F),\deg G)\sum_{p\ge P}\frac{\beta^{J+1}}{(1-\beta/P)(J+1)}\frac{1}{p^{J+1}}
    \\&\le
    \frac{6\max(\deg (G-F),\deg G)\beta^{J+1}}{(1-\beta/P)(J+1)}\biggl(
    \frac{1}{P^{J+1}}+\int_{P}^{\infty}\frac{dt}{t^{J+1}}\biggr)
    \\&\le
    \frac{6\max(\deg (G-F),\deg G)(\beta/P)^J\beta}{(1-\beta/P)(J+1)}\biggl(\frac{1}{P}+\frac{1}{J}\biggr),
\end{align*}
since $P\ge2$ and $J\ge3$.
We now approximate each sum over $p$ by using
Corollary~\ref{identityshort} and obtain
\begin{equation*}
    \sum_{\substack{p\ge P,\\ p\equiv a[q]}}
  \log\biggl(1-\frac{F(1/p)}{G(1/p)}\biggr)
  =
  \sum_{2\le j\le J}\bigl(b_{G-F}(j)-b_G(j)\bigr) Y_P(j;q,a|L)
  + I_3
\end{equation*}
where $I_3$ satisfies
\begin{align*}
    |I_3|
  &\le
   \sum_{2\le j\le J}|b_{G-F}(j)-b_G(j)|
  \frac{1}{P^{Lj}}
    +|I_2|
    \\&\le
    \sum_{2\le j\le
  J}4\max(\deg F,\deg G)\frac{\beta^j}{j}
  \frac{1}{P^{Lj}}
  +|I_2|.
\end{align*}
Therefore (and since $r\ge2$)
\begin{multline}
  \frac{|I_3|}{2\max(\deg F,\deg G)}\label{precbound}
  \le
  \frac{\beta^2(\beta/P)^{2L}}{1-\beta/P}
  +\frac{3(\beta/P)^J\beta}{(1-\beta/P)(J+1)}\biggl(\frac{1}{P}+\frac{1}{J}\biggr),
\end{multline}
and the choice $J=2L$ ends the proof.
\end{proof}

\section{Proof of Theorem~\ref{thm5}}

\begin{lem}
  \label{helper}
  We have $\binom{dN'}{dm'_1,\cdots,dm'_k}\ge
  \binom{N'}{m'_1,\cdots,m'_k}^d$.
\end{lem}

\begin{proof}
  The coefficient $\binom{dN'}{dm'_1,\cdots,dm'_k}$ is the number of
  partitions of a set of $dN'$ elements in parts of
  $dm'_1,\cdots,dm'_k$ elements. The product partitions are partitions.
\end{proof}

In \cite{Witt*37}, Witt proved a generalization of the Necklace
Identity which we present in the next lemma. 
\begin{lem}
  \label{MW}
  For $k\ge1$, we have
  \begin{equation}
    \label{Witt1937}
    1-
    \sum_{i=1}^k z_i=\prod_{\substack{m_1,\ldots,m_k\ge 0,\\ m_1+\ldots+m_k\ge1}}(1-z_1^{m_1}\cdots z_k^{m_k})^{M(m_1,\ldots,m_k)},
\end{equation}
where the integer $M(m_1,\ldots,m_k)$ is defined by
\begin{equation}
  \label{defM}
M(m_1,\ldots,m_k)=\frac{1}{N}
\sum_{d|\gcd(m_1,m_2,\ldots,m_k)}\mu(d)
\frac{(N/d)!}{(m_1/d)!\cdots (m_k/d)!}
\end{equation}
with $N=m_1+\ldots+m_k$. We have $M(m_1,\ldots,m_k)\le k^N/N$.
\end{lem}

\begin{proof}
  Only the bound needs to be proved as the identity may be found in
  \cite{Witt*37}. Each occuring multinomial
  is not more than $\binom{N}{m_1,\cdots,m_k}$ by
  Lemma~\ref{helper}. The multinomial Theorem concludes.
\end{proof}

\begin{proof}[Proof of Theorem~\ref{thm5}]
  Let $\Pi$ be the product to be computed.
  By employing Lemma~\ref{MW}, we
  find that
  \begin{equation*}
    1-\sum_{1\le \ell\le
      k}\frac{a_\ell}{p^{u_\ell s+v_\ell}}
    =\prod_{\substack{m_1,\ldots,m_k\ge 0,\\ m_1+\ldots+m_k\ge1}}
    \biggl(1-\frac{c(m_1,m_2,\ldots,m_k)}{p^{
        \sum_{\ell\le k}m_\ell (u_\ell s+v_\ell)}}\biggr)^{M(m_1,\ldots,m_k)},
  \end{equation*}
  with $ c(m_1,\ldots,m_k)$ given by
  \begin{equation}
    \label{defc}
    c(m_1,m_2,\ldots,m_k)
    =\prod_{\ell\le k}a_\ell^{m_\ell}.
  \end{equation}
  Each coefficient $ c(m_1,\ldots,m_k)$ is not more, in absolute value, than $A^N$, where
  $m_1+\ldots+m_k=N$.
  Note that, for each $\ell$, we have $\Re(u_\ell s+v_\ell)>1$, so
  that $\Re  \sum_{\ell\le k}m_\ell (u_\ell s+v_\ell)\ge
  m_1+\ldots+m_k=N$. It thus seems like a good idea to truncate the
  infinite product in~\eqref{infp} according to whether
  $m_1+\cdots+m_k=N\le N_0$ or not for some parameter $N_0\ge k$ that
  we will choose later.
  We readily find that, when $p\ge 2A$,
  \begin{multline*}
    \biggl|\log\prod_{\substack{m_1,\ldots,m_k\ge 0,\\ m_1+\ldots+m_k>N_0}}
    \biggl(1-\frac{c(m_1,m_2,\ldots,m_k)}{p^{
        \sum_{\ell\le k}m_\ell (u_\ell
        s+v_\ell)}}\biggr)^{M(m_1,\ldots,m_k)}\biggr|
    \\\le \frac32
    \sum_{\substack{m_1,\ldots,m_k\ge 0,\\ m_1+\ldots+m_k>N_0}}
    M(m_1,\ldots,m_k)\frac{A^N}{p^N}
    \\\le \frac32
    \sum_{\substack{N>N_0}}
    \binom{N+k}{k}\frac{(kA)^N}{N p^N}
  \end{multline*}
  as the number of solutions to $m_1+\ldots+m_k=N$ is the $N$-th
  coefficient of the power series $1/(1-z)^k$ which happens to be
  equal to $(1/k!)\frac{d}{dz^k}1/(1-z)$. We next check that, with
  $N=N_0+n+1$, we have
  $(n+1+N_0+k)\le (N_0+n+1)^2$ since $N_0\ge k$, and thus
  \begin{equation*}
    \frac{\binom{N+k}{k}}{N\binom{n+k}{k}}
    =\frac{(n+1+N_0+k)(n+N_0+k)\cdots(n+N_0+2)}{(n+k)(n+k-1)\cdots(n+1)\cdot(N_0+n+1)}
    \le \binom{N_0+k}{k}.
  \end{equation*}
  Hence, when $p\ge 2kA$, we
  have
  \begin{align*}
     \sum_{\substack{N>N_0}}
    \binom{N+k}{k}\frac{(kA)^N}{N p^N}
    &=
      \frac{(kA)^{N_0+1}}{p^{N_0+1}}
      \binom{N_0+k}{k}\sum_{n\ge
      0}
      \binom{n+k}{k}\frac{(kA)^n}{p^n}
    \\&\le
    \binom{N_0+k}{k}\frac{(kA)^{N_0+1}}{p^{N_0+1}}
    \frac{1}{(1-1/2)^k}.
  \end{align*}
  On summing over $p$, this yields
  \begin{equation}
    \label{infp}
    \Pi
    =I_1\prod_{\substack{m_1,\ldots,m_k\ge 0,\\ 1\le m_1+\ldots+m_k\le
      N_0}}
    \prod_{\substack{p\ge P,\\ p\equiv a[q]}}
    \biggl(1-\frac{c(m_1,m_2,\ldots,m_k)}{p^{
        \sum_{\ell\le k}m_\ell (u_\ell
        s+v_\ell)}}\biggr)^{M(m_1,\ldots,m_k)}
    ,
  \end{equation}
  where
  \begin{equation}
    \label{eq:2}
    |\log I_1|\le
    2^k\frac32\binom{N_0+k}{k}\frac{(kA)^{N_0+1}}{P^{N_0}}
    \biggl(\frac{1}{P}+\frac{1}{N_0}\biggr).
  \end{equation}
  We next note the following identity
  \begin{equation}
    \label{logc}
    \sum_{k\ge1}\frac{d^k}{kp^{kw}}
    =\sum_{f\ge1}\frac{\kappa_f(d)}{f}\sum_{{g\ge1}}\frac{1}{gp^{fgw}}
  \end{equation}
  where
  \begin{equation}
    \label{defkappa}
    \kappa_f(d)=
    \begin{cases}
      c&\quad\text{when $f=1$},\\
      c^{f}-c^{f-1}&\quad\text{when $f>1$}.
    \end{cases}
  \end{equation}
  We truncate identity~\eqref{logc} at $f\le F$ where $F$ is an integer, getting
  \begin{equation*}
    \label{eq:4}
     \sum_{k\ge1}\frac{d^k}{kp^{kw}}
    =\sum_{f\le
      F}\frac{\kappa_f(d)}{f}\sum_{{g\ge1}}\frac{1}{gp^{fgw}}
    +\Ocal^*\biggl(
    -\sum_{f>F}\frac{\max(1,|d|)^f}{f}
    \log\Bigl(1-p^{-f\Re w}\Bigr)
    \biggr).
  \end{equation*}
  We next use $-\log(1-x)\le 3x/2$ when $0\le x\le 1/2$. We assume
  that $p^{\Re w}\le 1/2$ and $p^{\Re w}\ge 2\max(1,|d|)$ to get
  \begin{equation*}
    -\sum_{f>F}\frac{\max(1,|d|)^f}{f}
    \log\Bigl(1-p^{-f\Re w}\Bigr)
    \le
    \frac32\sum_{f>F}\frac{\max(1,|d|)^f}{fp^{f\Re w}}
    \le \frac{3\max(1,|d|)^{F+1}}{(F+1)p^{(F+1)\Re w}}.
  \end{equation*}
  We have reached
  \begin{multline*}
    \prod_{\substack{p\ge P,\\ p\equiv a[q]}}
    \biggl(1-\frac{c(m_1,m_2,\ldots,m_k)}{p^{
        \sum_{\ell\le k}m_\ell (u_\ell
        s+v_\ell)}}\biggr)
    \\=
    \exp-\biggl\{\sum_{f\le
      F}\frac{\kappa_f(c(m_1,m_2,\ldots,m_k))}{f}
    \sum_{\substack{p\ge P,\\ p\equiv a[q]}}
    \log\Bigl(1-p^{-f\sum_{\ell\le k}m_\ell (u_\ell
      s+v_\ell)}\Bigr)
    \\+
        \Ocal^ *\biggl(
\frac{3\max(1,|c(m_1,m_2,\ldots,m_k)|)^{F+1}}{(F+1)P^{(F+1)\sum_{\ell\le k}m_\ell (u_\ell
      \sigma+v_\ell)}}\biggl(1+\frac{P}{F\sum_{\ell\le k}m_\ell (u_\ell
      \sigma+v_\ell)}\biggr)
    \biggr)
    \biggr\}
  \end{multline*}
  which simplifies info
  \begin{multline*}
    \prod_{\substack{p\ge P,\\ p\equiv a[q]}}
    \biggl(1-\frac{c(m_1,m_2,\ldots,m_k)}{p^{
        \sum_{\ell\le k}m_\ell (u_\ell
        s+v_\ell)}}\biggr)
    \\=
    \exp-\biggl\{\sum_{f\le
      F}\frac{\kappa_f(c(m_1,m_2,\ldots,m_k))}{f}
    \sum_{\substack{p\ge P,\\ p\equiv a[q]}}
    \log\Bigl(1-p^{-f\sum_{\ell\le k}m_\ell (u_\ell
      s+v_\ell)}\Bigr)
    \\+
    \frac{3A^{N(F+1)}}{(F+1)P^{(F+1)N}}\biggl(1+\frac{P}{FN}\biggr)
    \biggr\}.
  \end{multline*}
  We approximate the sum of the logs by~Corollary~\ref{identityshort}
  and get
  \begin{multline*}
    \prod_{\substack{p\ge P,\\ p\equiv a[q]}}
    \biggl(1-\frac{c(m_1,m_2,\ldots,m_k)}{p^{
        \sum_{\ell\le k}m_\ell (u_\ell
        s+v_\ell)}}\biggr)
    \\=
    \exp-\biggl\{\sum_{f\le
      F}
    \frac{\kappa_f(c(m_1,m_2,\ldots,m_k))}{f}
    Y_P\Bigl(\sum_{\ell\le k}m_\ell (u_\ell
    s+v_\ell);q,a|L\Bigr)
    \\+
    \Ocal^ *\biggl(
    \frac{A^{NF}(1+\log F)}{P^{LN}}
    +
    \frac{3A^{N(F+1)}}{(F+1)P^{(F+1)N}}\biggl(1+\frac{P}{FN}\biggr)
    \biggr)
    \biggr\}.
  \end{multline*}
  We then raise that to the power $M(m_1,m_2,\ldots,m_k)$ and sum over
  the $m_i$'s, getting, on recalling~\eqref{defZ},
  \begin{multline*}
    \Pi/I_1
    =
    \exp-Z
    \\+
    \Ocal^ *\biggl(
    \sum_{\substack{m_1,\ldots,m_k\ge 0,\\ 1\le m_1+\ldots+m_k\le
      N_0}}\frac{M(m_1,\ldots,m_k)A^{NF}(1+\log F)}{P^{LN}}
    \\+
    \sum_{\substack{m_1,\ldots,m_k\ge 0,\\ 1\le m_1+\ldots+m_k\le
      N_0}}\frac{3M(m_1,\ldots,m_k)A^{N(F+1)}}{(F+1)P^{(F+1)N}}\biggl(1+\frac{P}{FN}\biggr)
     \biggr)
   \biggr\}.
 \end{multline*}
 We now take $F=L$.
  The error term is bounded above by (since $P\ge 2kA$)
  \begin{equation*}
    \frac{kA^L}{P^L}\biggl(\frac{2^k}{k!}(1+\log
    L)+\frac{3\cdot2^kA}{k!(L+1)P}
    \biggl(1+\frac{P}{L}\biggr)\biggr).
  \end{equation*}
  We select $N_0=L$ and we gather our estimates to end the proof.
\end{proof}

\bibliographystyle{plain}

\end{document}